# A formula connecting the Bernoulli numbers with the Stieltjes constants

Donal F. Connon

dconnon@btopenworld.com

25 April 2011

**Abstract**

We show that the formula recently derived by Coffey for the Stieltjes constants in terms of the Bernoulli numbers is mathematically equivalent to the much earlier representation derived by Briggs and Chowla.

## 1. Introduction

The Stieltjes constants $\gamma_n(u)$ are the coefficients of the Laurent expansion of the Hurwitz zeta function $\varsigma(s,u)$ about $s=1$

$$(1.1) \qquad \varsigma(s,u) = \frac{1}{s-1} + \sum_{p=0}^{\infty} \frac{(-1)^p}{p!} \gamma_p(u)(s-1)^p$$

where $\gamma_p(u)$ are known as the generalised Stieltjes constants and we have [17]

$$(1.2) \qquad \gamma_0(u) = -\psi(u)$$

where $\psi(u)$ is the digamma function.

With $u=1$ in (1.1) the Hurwitz zeta function reduces to the Riemann zeta function

$$\varsigma(s) = \frac{1}{s-1} + \sum_{p=0}^{\infty} \frac{(-1)^p}{p!} \gamma_p (s-1)^p$$

The alternating Riemann zeta function $\varsigma_a(s)$ is defined for $\mathrm{Re}(s) > 0$ by

$$(1.3) \qquad \varsigma_a(s) = \sum_{n=1}^{\infty} \frac{(-1)^{n+1}}{n^s}$$

$$(1.4) \qquad = (1-2^{1-s})\varsigma(s)$$

We see from (1.3) that $\varsigma_a(1) = \log 2$ and from (1.4) we have

$$\log 2 = \lim_{s \to 1}(1 - 2^{1-s})\varsigma(s)$$

$$= \lim_{s \to 1} \frac{(1 - 2^{1-s})}{s - 1}(s - 1)\varsigma(s)$$

Using L'Hôpital's rule we find that $\lim_{s \to 1} \frac{(1 - 2^{1-s})}{s - 1} = \log 2$ and hence we easily determine the well-known limit

$$\lim_{s \to 1}(s - 1)\varsigma(s) = 1$$

We have the derivatives

$$\varsigma_a^{(k)}(1) = (-1)^k \sum_{n=1}^{\infty} \frac{(-1)^{n+1} \log^k n}{n^s}$$

## 2. A formula connecting the Bernoulli numbers with the Stieltjes constants

Noting that

$$(s - 1)\varsigma(s) = \left[\frac{s - 1}{1 - 2^{1-s}}\right][(1 - 2^{1-s})\varsigma(s)]$$

we see that

$$\left.\frac{d^{m+1}}{ds^{m+1}}[(s - 1)\varsigma(s)]\right|_{s=1} = \left.\frac{d^{m+1}}{ds^{m+1}}\left[\frac{s - 1}{1 - 2^{1-s}}\right][(1 - 2^{1-s})\varsigma(s)]\right|_{s=1}$$

$$= \left.\frac{d^{m+1}}{ds^{m+1}}\left[\frac{s - 1}{1 - 2^{1-s}}\right]\varsigma_a(s)]\right|_{s=1}$$

and applying the Leibniz rule for differentiation we obtain

$$= \sum_{l=0}^{m+1} \binom{m+1}{l} \left.\frac{d^{m+1-l}}{ds^{m+1-l}}\left[\frac{s - 1}{1 - 2^{1-s}}\right]\right|_{s=1} \varsigma_a^{(l)}(1)$$

We see from (1.1) that

$$\left.\frac{d^{m+1}}{ds^{m+1}}[(s - 1)\varsigma(s)]\right|_{s=1} = (-1)^m (m + 1)\gamma_m$$

and therefore we have

$$(2.1) \quad (-1)^m (m+1)\gamma_m = \sum_{l=0}^{m+1} \binom{m+1}{l} \frac{d^{m+1-l}}{ds^{m+1-l}} \left[ \frac{s-1}{1-2^{1-s}} \right]_{s=1} \varsigma_a^{(l)}(1)$$

The Bernoulli numbers $B_n$ are given by the generating function [9]

$$\frac{t}{e^t - 1} = \sum_{n=0}^{\infty} B_n \frac{t^n}{n!} \quad , \quad (|t| < 2\pi)$$

and we therefore have

$$\frac{1}{x^{1-s} - 1} = \frac{1}{\exp[-(s-1)\log x] - 1} = -\frac{1}{(s-1)\log x} + \sum_{k=0}^{\infty} \frac{(-1)^k B_{k+1} \log^k x}{(k+1)!} (s-1)^k$$

which is valid for $|s-1| < \dfrac{2\pi}{\log x}$.

Hence we have

$$\frac{s-1}{1-x^{1-s}} = \frac{1}{\log x} + \sum_{k=0}^{\infty} \frac{(-1)^{k+1} B_{k+1} \log^k x}{(k+1)!} (s-1)^{k+1}$$

and differentiation results in

$$\frac{d^{m+1-l}}{ds^{m+1-l}} \left[ \frac{s-1}{1-x^{1-s}} \right] = \sum_{k=0}^{\infty} \frac{(-1)^{k+1} B_{k+1} \log^k x}{(k+1)!} (k+1)k \cdots (k-m+l+1)(s-1)^{k-m+l}$$

so that

$$\frac{d^{m+1-l}}{ds^{m+1-l}} \left[ \frac{s-1}{1-x^{1-s}} \right]_{s=1} = (-1)^{m-l+1} B_{m-l+1} \log^{m-l} x$$

and in particular we have

$$\frac{d^{m+1-l}}{ds^{m+1-l}} \left[ \frac{s-1}{1-2^{1-s}} \right]_{s=1} = (-1)^{m-l+1} B_{m-l+1} \log^{m-l} 2$$

We then see from (2.1) that

$$(2.2) \quad (m+1)\gamma_m = -\sum_{l=0}^{m+1} \binom{m+1}{l} (-1)^l B_{m-l+1} \varsigma_a^{(l)}(1) \log^{m-l} 2$$

In 2006 Coffey ([3] and [5]) showed that

$$(2.3) \quad \gamma_m = -m! \sum_{l=1}^{m+1} \frac{B_{m-l+1}}{(m-l+1)!} \frac{\log^{m-l} 2}{l!} (-1)^l \varsigma_a^{(l)}(1) - \frac{B_{m+1}}{m+1} \log^{m+1} 2$$

and using some basic binomial number identities we see that

$$= -\sum_{l=1}^{m+1} \binom{m}{l} \frac{B_{m-l+1}}{m-l+1} \log^{m-l} 2 (-1)^l \varsigma_a^{(l)}(1) - \frac{B_{m+1}}{m+1} \log^{m+1} 2$$

and we therefore obtain

$$\gamma_m = -\frac{1}{m+1} \sum_{l=1}^{m+1} \binom{m+1}{l} B_{m-l+1} \log^{m-l} 2 (-1)^l \varsigma_a^{(l)}(1) - \frac{B_{m+1}}{m+1} \log^{m+1} 2$$

Since $\varsigma_a^{(0)}(1) = \log 2$, with the summation starting at $l = 0$, we may write this as

$$(2.4) \quad \gamma_m = -\frac{1}{m+1} \sum_{l=0}^{m+1} \binom{m+1}{l} B_{m-l+1} \log^{m-l} 2 \cdot (-1)^l \varsigma_a^{(l)}(1)$$

and this is the formula reported by Zhang and Williams [17] in 1994 (and this corresponds with (2.2)).

Reindexing (2.4) so that $l = m+1-k$ gives us

$$(2.4.1) \quad \gamma_m = -\frac{1}{m+1} \sum_{k=0}^{m+1} \binom{m+1}{k} B_k \log^{k-1} 2 \cdot (-1)^{m+1-k} \varsigma_a^{(m+1-k)}(1)$$

and this representation is employed later to prove Kluyver's formula (5.1).

The expression originally derived by Liang and Todd [12] in 1972 was

$$(2.5) \quad \gamma_m = \frac{\log^{m+1} 2}{m+1} \sum_{l=1}^{m+1} \binom{m+1}{l} \frac{B_{m-l+1}}{l+1} + \frac{1}{m+1} \sum_{l=1}^{m+1} \binom{m+1}{l} B_{m-l+1} \log^{m-l} 2 \cdot (-1)^{l+1} \varsigma_a^{(l)}(1)$$

which, at first glance, appears very different from (2.4) but, as will be seen below, they are in fact the same.

Zhang and Williams [17] noted that for $m = 1, 2, ...$ we have

$$(2.6) \quad \sum_{l=1}^{m} \binom{m}{l} \frac{B_{m-l}}{l+1} = 0$$

and a slightly modified version of their proof is shown below.

We recall that the Bernoulli polynomials may be expressed by [9, p.2]

(2.7) $$B_m(x) = \sum_{l=0}^{m} \binom{m}{l} B_{m-l} x^l$$

where integration results in

$$\int_0^1 B_m(x)dx = \sum_{l=0}^{m} \binom{m}{l} \frac{B_{m-l}}{l+1}$$

Alternatively we have

$$\int_0^1 B_m(x)dx = \frac{1}{m+1}[B_{m+1}(1) - B_{m+1}(0)] = 0$$

and therefore we easily deduce (2.6).

We accordingly see that

$$\sum_{l=1}^{m+1} \binom{m+1}{l} \frac{B_{m-l+1}}{l+1} = \sum_{l=0}^{m+1} \binom{m+1}{l} \frac{B_{m-l+1}}{l+1} - B_{m+1} = -B_{m+1}$$

and hence we may write the Liang and Todd equation (2.5) as

$$\gamma_m = -\frac{B_{m+1} \log^{m+1} 2}{m+1} + \frac{1}{m+1} \sum_{l=1}^{m+1} \binom{m+1}{l} B_{m-l+1} \log^{m-l} 2 \cdot (-1)^{l+1} \varsigma_a^{(l)}(1)$$

$$= -\frac{1}{m+1} \sum_{l=0}^{m+1} \binom{m+1}{l} B_{m-l+1} \log^{m-l} 2 \cdot (-1)^l \varsigma_a^{(l)}(1)$$

where the summation starts at $l = 0$ and we note again that $\varsigma_a^{(0)}(1) = \varsigma_a(1) = \log 2$. Therefore, we have shown that (2.4) and (2.5) are equivalent.

## 3. Equivalence of the Briggs and Chowla formula with the Coffey representation

The following result was obtained by Briggs and Chowla [2] in 1955

(3.1) $$\varsigma_a^{(k)}(1) = k! \sum_{r=1}^{k+1} \frac{(-1)^{r+1} \log^r 2}{r!} A_{k-r}$$

with $A_n = \dfrac{(-1)^n}{n!}\gamma_n$ and $A_{-1} = 1$.

and, as shown by Dilcher [8], this may also be expressed as

(3.2) $\quad (-1)^l \varsigma_a^{(l)}(1) = \dfrac{1}{l+1}\log^{l+1} 2 - \displaystyle\sum_{k=0}^{l-1}\binom{l}{k}\gamma_k \log^{l-k} 2$

Other proofs have been provided, inter alia, by Zhang [16] and the author [7]. Examples of (3.2) are set out below:

(3.3) $\quad -\varsigma_a^{(1)}(1) = \displaystyle\sum_{k=1}^{\infty}(-1)^{k-1}\dfrac{\log k}{k} = \dfrac{1}{2}\log^2 2 - \gamma \log 2$

(3.4) $\quad \varsigma_a^{(2)}(1) = \displaystyle\sum_{k=1}^{\infty}(-1)^{k-1}\dfrac{\log^2 k}{k} = \dfrac{1}{3}\log^3 2 - \gamma \log^2 2 - 2\gamma_1 \log 2$

Equation (3.3) was posed as a problem by Klamkin [10] in 1954 and is closely related to an earlier problem posed by Sandham [15] in 1950.

In 2009, to the author's surprise [7], it was noted that substituting (3.2) in (2.4) did not appear to produce any additional information

$$\gamma_m = -\dfrac{1}{m+1}\sum_{l=0}^{m+1}\binom{m+1}{l}B_{m-l+1}\log^{m-l} 2 \cdot \left[\dfrac{1}{l+1}\log^{l+1} 2 - \sum_{k=0}^{l-1}\binom{l}{k}\gamma_k \log^{l-k} 2\right]$$

For example, with $m = 1$ we simply find that $\gamma_1 = \gamma_1$! As shown below, the reason for this is that (2.4) and (3.2) are in fact equivalent representations of the same thing.

In 1964 Riordan [13] reported the following inverse relations involving the Bernoulli numbers

(3.5) $\quad a_n = \displaystyle\sum_{k=0}^{n}\binom{n}{k}\dfrac{b_k}{n-k+1} \Leftrightarrow b_n = \displaystyle\sum_{k=0}^{n}\binom{n}{k}B_{n-k}a_k$

and we note that a simple reindexing shows that

(3.6) $\quad b_n = \displaystyle\sum_{k=0}^{n}\binom{n}{k}B_{n-k}a_k = \displaystyle\sum_{k=0}^{n}\binom{n}{k}B_k a_{n-k}$

Using (2.4) we have

$$b_n \equiv \frac{-n\gamma_{n-1}}{\log^{n-1} 2} = \sum_{k=0}^{n} \binom{n}{k} B_{n-k} \frac{(-1)^k \varsigma_a^{(k)}(1)}{\log^k 2}$$

and inverting this with Riordan's formula therefore gives us

$$\frac{(-1)^n \varsigma_a^{(n)}(1)}{\log^n 2} = -\sum_{k=0}^{n} \binom{n}{k} \frac{k\gamma_{k-1}}{(n-k+1)\log^{k-1} 2}$$

so that we have

$$(-1)^n \varsigma_a^{(n)}(1) = -\sum_{k=0}^{n} \binom{n}{k} \frac{k\gamma_{k-1} \log^{n-k+1} 2}{n-k+1}$$

As we shall see below, this is the same as the formula (3.2) derived by Briggs and Chowla.

Reindexing $k \to n-r$ gives us

$$(-1)^n \varsigma_a^{(n)}(1) = -\sum_{r=0}^{n} \binom{n}{n-r} \frac{(n-r)\gamma_{n-r-1} \log^{r+1} 2}{r+1}$$

$$= -\sum_{r=0}^{n} \binom{n}{r} \frac{(n-r)\gamma_{n-r-1} \log^{r+1} 2}{r+1}$$

A further reindexing $r \to p-1$ results in

$$= -\sum_{p=1}^{n+1} \binom{n}{p-1} \frac{(n-p)\gamma_{n-p} \log^p 2}{p}$$

$$= \frac{\log^{n+1} 2}{n+1} - \sum_{p=1}^{n} \binom{n}{p-1} \frac{(n-p)\gamma_{n-p} \log^p 2}{p}$$

and isolating the term $p = n$ gives us

$$= \frac{\log^{n+1} 2}{n+1} - \sum_{p=1}^{n-1} \binom{n}{p-1} \frac{(n-p)\gamma_{n-p} \log^p 2}{p}$$

Since $\binom{n}{-1} = 0$ we have

$$= \frac{\log^{n+1} 2}{n+1} - \sum_{p=0}^{n-1} \binom{n}{p-1} \frac{(n-p)\gamma_{n-p} \log^p 2}{p}$$

We have by reindexing $p \to n-q$

$$= \frac{\log^{n+1} 2}{n+1} - \sum_{q=0}^{n-1} \binom{n}{n-q-1} \frac{q\gamma_q \log^{n-q} 2}{n-q}$$

and finally we obtain the Briggs and Chowla formula

$$(-1)^n \varsigma_a^{(n)}(1) = \frac{\log^{n+1} 2}{n+1} - \sum_{q=0}^{n-1} \binom{n}{q} \gamma_q \log^{n-q} 2$$

□

Using a series acceleration technique originally devised by Amore [1], Coffey [4] also showed that for integers $m \geq 1$ and Re $\lambda > 0$

$$\gamma_m = -\frac{m!}{1+\lambda} \sum_{n=1}^{m} \frac{B_{m-n+1}}{(m-n+1)!} \frac{\log^{m-n} 2}{n!} \sum_{k=0}^{\infty} \left(\frac{\lambda}{1+\lambda}\right)^k \sum_{j=0}^{k} \binom{k}{j} (-1)^j \frac{1}{\lambda^j} \frac{\log^n(j+1)}{j+1}$$

$$- \frac{1}{(1+\lambda)\log 2} \frac{1}{(m+1)} \sum_{k=0}^{\infty} \left(\frac{\lambda}{1+\lambda}\right)^k \sum_{j=0}^{k} \binom{k}{j} (-1)^j \frac{\log^{m+1}(j+1)}{j+1} - \frac{B_{m+1} \log^{m+1} 2}{m+1}$$

With $\lambda = 1$ we have

$$\gamma_m = -m! \sum_{n=1}^{m} \frac{B_{m-n+1}}{(m-n+1)!} \frac{\log^{m-n} 2}{n!} \sum_{k=0}^{\infty} \frac{1}{2^{k+1}} \sum_{j=0}^{k} \binom{k}{j} (-1)^j \frac{\log^n(j+1)}{j+1}$$

$$- \frac{1}{(m+1)\log 2} \sum_{k=0}^{\infty} \frac{1}{2^{k+1}} \sum_{j=0}^{k} \binom{k}{j} (-1)^j \frac{\log^{m+1}(j+1)}{j+1} - \frac{B_{m+1} \log^{m+1} 2}{m+1}$$

and using the notation

$$a_n = \log^{-n} 2 \sum_{k=0}^{\infty} \frac{1}{2^{k+1}} \sum_{j=0}^{k} \binom{k}{j} (-1)^j \frac{\log^n(j+1)}{j+1}$$

we have

$$\gamma_m = -\frac{\log^m 2}{m+1} \sum_{n=1}^{m} \binom{m+1}{n} \frac{B_{m-n+1} a_n}{m-n+1} - \frac{a_{m+1} \log^m 2}{m+1} - \frac{B_{m+1} \log^{m+1} 2}{m+1}$$

$$\text{(3.7)} \qquad = -\frac{\log^m 2}{m+1} \sum_{n=0}^{m+1} \binom{m+1}{n} \frac{B_{m-n+1} a_n}{m-n+1}$$

We now refer to the Hasse identity for the alternating Hurwitz zeta function (see equation (4.4.79) in [6])

$$\text{(3.8)} \qquad \varsigma_a(s,x) = \sum_{k=0}^{\infty} \frac{1}{2^{k+1}} \sum_{j=0}^{k} \binom{k}{j} \frac{(-1)^j}{(x+j)^s}$$

Differentiation gives us

$$\frac{\partial^n}{\partial s^n} \varsigma_a(s,x) = (-1)^n \sum_{k=0}^{\infty} \frac{1}{2^{k+1}} \sum_{j=0}^{k} \binom{k}{j} \frac{(-1)^j \log^n(x+j)}{(x+j)^s}$$

and we see that

$$a_n = (-1)^n \varsigma_a^{(n)}(1) \log^{-n} 2$$

$$a_0 = \varsigma_a^{(0)}(1) = \log 2$$

Substituting this in (3.7) results in

$$\text{(3.9)} \qquad \gamma_m = -\frac{1}{m+1} \sum_{n=0}^{m+1} \binom{m+1}{n} \frac{B_{m-n+1}(-1)^n \varsigma_a^{(n)}(1) \log^{m-n} 2}{m-n+1}$$

which concurs with (2.4).

**5. Kluyver's formula for the Stieltjes constants**

Kluyver [11] showed in 1927 that the Stieltjes constants may be represented in terms of the Bernoulli polynomials as follows

$$\text{(5.1)} \qquad \gamma_p = \frac{\log^p 2}{p+1} \sum_{k=1}^{\infty} \frac{(-1)^k}{k} B_{p+1}\left(\frac{\log k}{\log 2}\right)$$

We may obtain this formula in the following manner. We denote $c_p$ as

$$c_p = \frac{\log^p 2}{p+1} \sum_{k=1}^{\infty} \frac{(-1)^k}{k} B_{p+1}\left(\frac{\log k}{\log 2}\right)$$

and using (2.7) we have

$$B_{p+1}(x) = \sum_{r=0}^{p+1} \binom{p+1}{r} B_r x^{p+1-r}$$

and we obtain

$$c_p = \frac{\log^p 2}{p+1} \sum_{k=1}^{\infty} \frac{(-1)^k}{k} \sum_{r=0}^{p+1} \binom{p+1}{r} B_r \left[\frac{\log k}{\log 2}\right]^{p+1-r}$$

where we have employed square brackets to make it clear that $c_p$ is now being expressed in terms of the Bernoulli numbers. We see that

$$c_p = \frac{\log^p 2}{p+1} \sum_{r=0}^{p+1} \binom{p+1}{r} B_r \sum_{k=1}^{\infty} \frac{(-1)^k}{k} \left[\frac{\log k}{\log 2}\right]^{p+1-r}$$

$$= \frac{1}{p+1} \sum_{r=0}^{p+1} \binom{p+1}{r} B_r \log^{r-1} 2 \sum_{k=1}^{\infty} \frac{(-1)^k}{k} \log^{p+1-r} k$$

and substituting

$$\varsigma_a^{(j)}(1) = (-1)^j \sum_{k=1}^{\infty} \frac{(-1)^{k+1} \log^j k}{k^s}$$

we obtain

$$c_p = -\frac{1}{p+1} \sum_{r=0}^{p+1} \binom{p+1}{r} B_r \log^{r-1} 2 (-1)^{p+1-r} \varsigma_a^{(p+1-r)}(1)$$

Comparing this with (2.4.1) we see that $c_p = \gamma_p$ and we have therefore proved Kluyver's formula (5.1).

□

We may extend Kluyver's formula (5.1) in the following manner. We denote $c_p(u)$ as

$$c_p(u) = \frac{\log^p 2}{p+1} \sum_{k=1}^{\infty} \frac{(-1)^k}{k} B_{p+1}\left(u \frac{\log k}{\log 2}\right)$$

and using (2.7) we have

$$c_p(u) = \frac{\log^p 2}{p+1} \sum_{k=1}^{\infty} \frac{(-1)^k}{k} \sum_{r=0}^{p+1} \binom{p+1}{r} B_r \left[u \frac{\log k}{\log 2}\right]^{p+1-r}$$

We see that

$$c_p(u) = \frac{\log^p 2}{p+1} \sum_{r=0}^{p+1} \binom{p+1}{r} B_r \sum_{k=1}^{\infty} \frac{(-1)^k}{k} \left[ u \frac{\log k}{\log 2} \right]^{p+1-r}$$

$$= \frac{\log^p 2}{p+1} \sum_{r=0}^{p+1} \binom{p+1}{r} B_r \frac{u^{p+1-r}}{\log^{p+1-r} 2} \sum_{k=1}^{\infty} \frac{(-1)^k}{k} \log^{p+1-r} k$$

and substituting

$$\varsigma_a^{(j)}(1) = (-1)^j \sum_{k=1}^{\infty} \frac{(-1)^{k+1} \log^j k}{k^s}$$

we obtain

$$c_p(u) = -\frac{1}{p+1} \sum_{r=0}^{p+1} \binom{p+1}{r} B_r \log^{r-1} 2 (-1)^{p+1-r} u^{p+1-r} \varsigma_a^{(p+1-r)}(1)$$

$\square$

Before I accidentally, but rather fortuitously, came across Riordan's inverse relation (3.5), I tried to analyse Kluyver's formula as follows:

Starting with

$$\gamma_p = \frac{\log^p 2}{p+1} \sum_{k=1}^{\infty} \frac{(-1)^k}{k} B_{p+1} \left( \frac{\log k}{\log 2} \right)$$

and substituting the expression for the Bernoulli polynomials

$$B_{p+1}(x) = \sum_{r=0}^{p+1} \binom{p+1}{r} B_r x^{p+1-r}$$

we have as before

$$\gamma_p = \frac{\log^p 2}{p+1} \sum_{k=1}^{\infty} \frac{(-1)^k}{k} \sum_{r=0}^{p+1} \binom{p+1}{r} B_r \left[ \frac{\log k}{\log 2} \right]^{p+1-r}$$

which is now expressed in terms of the Bernoulli numbers

$$= \frac{\log^p 2}{p+1} \sum_{r=0}^{p+1} \binom{p+1}{r} B_r \sum_{k=1}^{\infty} \frac{(-1)^k}{k} \left[ \frac{\log k}{\log 2} \right]^{p+1-r}$$

$$= \frac{\log^p 2}{p+1} \left( \frac{1}{\log 2} \right)^{p+1} \sum_{r=0}^{p+1} \binom{p+1}{r} \left( \frac{1}{\log 2} \right)^{-r} B_r \sum_{k=1}^{\infty} \frac{(-1)^k}{k} \log^{p+1-r} k$$

$$= \frac{1}{(p+1)\log 2} \sum_{r=0}^{p+1} \binom{p+1}{r} \log^r 2 \, B_r \sum_{k=1}^{\infty} \frac{(-1)^k}{k} \log^{p+1-r} k$$

$$= \frac{1}{(p+1)\log 2} \sum_{r=0}^{p} \binom{p+1}{r} \log^r 2 \, B_r \sum_{k=1}^{\infty} \frac{(-1)^k}{k} \log^{p+1-r} k$$

$$+ \frac{1}{(p+1)\log 2} \log^{p+1} 2 \, B_{p+1} \sum_{k=1}^{\infty} \frac{(-1)^k}{k}$$

$$= \frac{1}{(p+1)\log 2} \sum_{r=0}^{p} \binom{p+1}{r} \log^r 2 \, B_r \sum_{k=1}^{\infty} \frac{(-1)^k}{k} \log^{p+1-r} k$$

$$- \frac{\log^{p+1} 2}{p+1} B_{p+1}$$

Then substituting the Briggs and Chowla formula (3.2)

$$\sum_{k=1}^{\infty} \frac{(-1)^k}{k} \log^n k = \sum_{m=0}^{n-1} \binom{n}{m} \gamma_m \log^{n-m} 2 - \frac{\log^{n+1} 2}{n+1}$$

gives us

$$(p+1)\gamma_p = \sum_{r=0}^{p} \binom{p+1}{r} B_r \sum_{m=0}^{p-r} \binom{p+1-r}{m} \gamma_m \log^{p-m} 2 - \log^{p+1} 2 \sum_{r=0}^{p} \binom{p+1}{r} \frac{B_r}{p-r+2} - \log^{p+1} 2 \, B_{p+1}$$

We write

$$\sum_{r=0}^{p} \binom{p+1}{r} \frac{B_r}{p-r+2} = \sum_{r=0}^{p+1} \binom{p+1}{r} \frac{B_r}{p-r+2} - B_{p+1}$$

and noting that $B_{p+1}(x) = \sum_{r=0}^{p+1} \binom{p+1}{r} B_r x^{p+1-r}$ we obtain

$$\int_0^1 B_{p+1}(x) dx = \sum_{r=0}^{p+1} \binom{p+1}{r} \frac{B_r}{p-r+2}$$

Using

$$\int_0^1 B_{p+1}(x) dx = B_{p+2}(1) - B_{p+2}(0) = 0$$

we therefore see that

$$\sum_{r=0}^{p+1}\binom{p+1}{r}\frac{B_r}{p-r+2}=0$$

Hence we obtain

$$(p+1)\gamma_p = \sum_{r=0}^{p}\binom{p+1}{r}B_r\sum_{m=0}^{p-r}\binom{p+1-r}{m}\gamma_m \log^{p-m} 2$$

which turned out to be a rather useless identity because, for example, letting $p=1$ we simply end up with $2\gamma_1 = 2\gamma_1$.

## 6. Some applications of Riordan's inversion formulae

As a minor diversion, we consider other examples of Riordan's inversion formula (3.5)

$$a_n = \sum_{k=0}^{n}\binom{n}{k}\frac{b_k}{n-k+1} \Leftrightarrow b_n = \sum_{k=0}^{n}\binom{n}{k}B_{n-k}a_k$$

in the case where $a_k = (-1)^k B_k$ so that

$$b_n = \sum_{k=0}^{n}\binom{n}{k}B_{n-k}(-1)^k B_k$$

We see that

$$\sum_{k=0}^{n}\binom{n}{k}B_{n-k}(-1)^k B_k = B_0 B_n + \sum_{k=1}^{n}\binom{n}{k}B_{n-k}(-1)^k B_k$$

$$= B_n + \sum_{k=1}^{n}\binom{n}{k}B_{n-k}(-1)^k B_k$$

and noting an identity reported by Rubenstein [14]

$$\sum_{k=1}^{n}\binom{n}{k}B_{n-k}(-1)^k B_k = -nB_n$$

we obtain

$$b_n = \sum_{k=0}^{n}\binom{n}{k}B_{n-k}(-1)^k B_k = (1-n)B_n$$

Then using (3.5) we end up with a combinatorial identity

$$(-1)^n B_n = \sum_{k=0}^{n} \binom{n}{k} \frac{(1-k)B_k}{n-k+1}$$

$$= \sum_{k=0}^{n} \binom{n}{k} B_k - n \sum_{k=0}^{n} \binom{n}{k} \frac{B_k}{n-k+1}$$

$$= \sum_{k=0}^{n} \binom{n}{k} B_k - n \sum_{k=0}^{n-1} \binom{n}{k} \frac{B_k}{n-k+1} - nB_n$$

We have the identity [9] for $p \geq 1$

$$\sum_{k=0}^{p-1} \binom{p}{k} \frac{B_k}{p-k+1} = -B_p$$

and for $p = n$ this becomes for $n \geq 1$

$$\sum_{k=0}^{n-1} \binom{n}{k} \frac{B_k}{n-k+1} = -B_n$$

Therefore we have

$$(-1)^n B_n = \sum_{k=0}^{n} \binom{n}{k} B_k$$

and using

$$B_n = \sum_{k=0}^{n} \binom{n}{k} B_k$$

we obtain the rather obvious result

$$(-1)^n B_n = B_n$$

## 7. Another identity involving the Stieltjes constants

In a paper recently submitted to arXiv, "Some integrals involving the Stieltjes constants: Part II", we showed that

(7.1) $$\varsigma_a^{(l)}(1,x) = \frac{1}{2} \sum_{k=0}^{l} \binom{l}{k} \log^{n-k} 2 \cdot \left[ \gamma_k\left(\frac{x}{2}\right) - \gamma_k\left(\frac{1+x}{2}\right) \right]$$

so that

(7.2) $$\varsigma_a^{(l)}(1) = \frac{1}{2}\sum_{k=0}^{l}\binom{l}{k}\log^{l-k}2\cdot\left[\gamma_k\left(\frac{1}{2}\right)-\gamma_k\right]$$

Substituting this in (3.2)

$$(-1)^l \varsigma_a^{(l)}(1) = \frac{1}{l+1}\log^{l+1}2 - \sum_{k=0}^{l-1}\binom{l}{k}\gamma_k \log^{l-k}2$$

results in

(7.3) $$(-1)^l \frac{1}{2}\sum_{k=0}^{l}\binom{l}{k}\log^{-k}2\cdot\left[\gamma_k\left(\frac{1}{2}\right)-\gamma_k\right] = \frac{1}{l+1}\log 2 - \sum_{k=0}^{l}\binom{l}{k}\gamma_k \log^{-k}2 + \gamma_l \log^{-l}2$$

and, as expected, with $l = 0$ we get the well-known value $\psi\left(\frac{1}{2}\right) = -\gamma - 2\log 2$.

Then substituting the well-known formula for $\gamma_k\left(\frac{1}{2}\right)$

(7.4) $$\gamma_k\left(\frac{1}{2}\right) = -\gamma_k + 2(-1)^k \frac{\log^{k+1}2}{k+1} + 2\sum_{j=0}^{k}\binom{k}{j}(-1)^j \gamma_{k-j}\log^j 2$$

gives us

$$(-1)^l \log 2 \sum_{k=0}^{l}\binom{l}{k}\frac{(-1)^k}{k+1} + (-1)^l \sum_{k=0}^{l}\binom{l}{k}\log^{-k}2\sum_{j=0}^{k}\binom{k}{j}(-1)^j \gamma_{k-j}\log^j 2$$

$$+\left[1-(-1)^l\right]\sum_{k=0}^{l}\binom{l}{k}\gamma_k \log^{-k}2 = \frac{1}{l+1}\log 2 + \gamma_l \log^{-l}2$$

and since $\sum_{k=0}^{l}\binom{l}{k}\frac{(-1)^k}{k+1} = \frac{1}{l+1}$ we obtain the curious identity

(7.5) $$(-1)^l \sum_{k=0}^{l}\binom{l}{k}\log^{-k}2\sum_{j=0}^{k}\binom{k}{j}(-1)^j \gamma_{k-j}\log^j 2 + \left[1-(-1)^l\right]\sum_{k=0}^{l}\binom{l}{k}\gamma_k \log^{-k}2$$

$$= \frac{\left[1-(-1)^l\right]}{l+1}\log 2 + \gamma_l \log^{-l}2$$

Donal F. Connon
Elmhurst
Dundle Road
Matfield
Kent TN12 7HD
dconnon@btopenworld.com